\documentclass[12pt]{amsart}
\usepackage{amssymb}
\usepackage{enumitem}
\usepackage{hyperref}
\usepackage{tabularx}
\usepackage{booktabs}
\usepackage{caption}
\usepackage[T2A,T1]{fontenc}
\usepackage[utf8]{inputenc}
\usepackage[russian,english]{babel}

\newtheorem*{prob*}{Open problem}

\theoremstyle{definition}

\theoremstyle{remark}

\newtheorem*{rem*}{Remark}

\newcommand{\kringel}{\mathbin{\raise1pt\hbox{$\scriptstyle\circ$}}}
\newcommand{\pkt}{\mathbin{\raise0pt\hbox{$\scriptstyle\bullet$}}}

\renewcommand{\phi}{\varphi}

\begin{document}


\title[Igor Dmitrievich Ado]{In memory of Igor Dmitrievich Ado}

\author[D. Burde]{Dietrich Burde}
\author[V. Gubarev]{Vsevolod Gubarev}
\address{Fakult\"at f\"ur Mathematik\\
Universit\"at Wien\\
  Oskar-Morgenstern-Platz 1\\
  1090 Wien \\
  Austria}
\email{dietrich.burde@univie.ac.at}
\address{Fakult\"at f\"ur Mathematik\\
Universit\"at Wien\\
  Oskar-Morgenstern-Platz 1\\
  1090 Wien \\
  Austria}
\email{vsevolod.gubarev@univie.ac.at, wsewolod89@gmail.com}

\date{\today}


\begin{abstract}
We give a translation from Russian into English of the article 
\begin{otherlanguage*}{russian}
"Памяти Игоря Дмитриевича Адо"
\end{otherlanguage*}
written by A.V. Dorodnov and I.I. Sakhaev and published
in Izv.\ Vyssh.\ Uchebn.\ Zaved.\ Mat.\ no. 8, (1984), 87--88.
It is an orbituary for I. D. Ado. A translation might be useful in general, and in particular
for a possible Wikipedia entry of Ado's life in English. In the references we list all known 
$11$ publications of I.D. Ado, taken from the article and the MATHSCINET of the AMS. 
The original orbituary only lists $9$ publications. 
\end{abstract}

\maketitle

On June $29$th $1983$ the famous Soviet mathematician, doctor of physical and mathematical 
sciences and Professor Igor Dmitrievich Ado passed away at the age of $73$.
I.D. Ado was born in Kazan in January $1910$ into the family of a state employee and he lived 
in Kazan till the end of his life. After leaving school Igor Dmitrievich entered the 
faculty of mathematics and physics at Kazan State University, named after 
V.I. Lenin, from which he graduated successfully in $1931$. He was admitted to the PhD 
study at the Chair of Mathematics (since 1934 – Chair of Algebra) under the supervision of N.G. Chebotarev. 
Igor Dmitrievich finished successfully his PhD study by preparing a scientific qualifying work for
the degree of a Candidate (PhD) of physical and mathematical sciences. 
The University board awarded him for this work the degree of {\em Doctor nauk} 
(doctor of sciences) of physic-mathematical sciences. Igor Dmitrievich solved in his thesis a 
current problem of modern 
algebra connected to representation theory of Lie algebras and Lie groups. 
More precisely, he obtained the result which is now known as {\em Ado's Theorem}: every finite-dimensional 
Lie algebra over a field of characteristic zero has a faithful finite-dimensional linear 
representation. \\[0.3cm]
In $1932$, during the International Mathematical Congress in Z\"urich, N.G. Chebotarev realized in 
the conversation with van der Waerden that the problem of representations of finite-dimensional 
Lie algebras was still open. After that N.G. Chebotarev suggested this problem to I.D. Ado as the 
theme of his future thesis. Igor Dmitrievich solved the problem and obtained a brilliant result, which brought 
him worldwide fame. \\[0.3cm]
Ado' Theorem has attracted the attention of many famous mathematicians who tried to improve its proof. 
In $1938$ \'E. Cartan gave another proof and in $1947$ I.D. Ado found a new proof of his theorem \cite{ADO7}. 
In $1937$ G. Birkhoff proved Ado's Theorem for nilpotent Lie algebras. In $1948$ K. Iwasawa proved the 
theorem for a field of positive characteristic $p>0$, so that the result about faithful linear representation 
of finite-dimensional Lie algebras is now also called the {\em Ado-Iwasawa Theorem}. \\[0.5cm]
I.D. Ado got a series of results on the structure of finite continuous groups \cite{ADO1}, 
on representations of finite continuous groups by linear substitutions \cite{ADO2} and
on nilpotent algebras and $p$-groups \cite{ADO3}. He had a significant interest in group theory and stated some 
results in character theory of finite groups \cite{ADO4} and in local finite $p$-groups with minimality condition 
for normal factors \cite{ADO5,ADO6,ADO8}. His last scientific work \cite{ADO9} without coauthors was devoted to 
the theory of linear representations of finite groups. \\[0.3cm]
After the defense of his doctoral thesis Igor Dmitrievich started to work at Kazan State University. 
Since $1936$ till $1942$ he held the position of a professor at the Chair of Algebra. 
In $1942$ I.D. Ado moved to the Kazan State Chemical Technological Institute named after S.M. Kirov, 
where he hold the position of a the Chair of High Mathematics until his death. 
He hold the position of a professor form $1942$ to $1958$ and from $1970$ to $1983$. During the period from 
$1958$ to $1970$ he was the head of the chair. \\[0.3cm]
Igor Dmitrievich was a wonderful teacher. He gave lectures on high theoretical level and at the same 
time comprehensible for students. On seminars he solved the tasks pretending that he saw them for 
the first time. He claimed hypotheses, checked them, highlighted particular cases and led 
students to the final solution. He thought out loud and taught students to think. 
Students and colleagues respected him and loved him.

\section*{Acknowledgments}
Dietrich Burde is supported by the Austrian Science Foun\-da\-tion FWF, grant I3248.
Vsevolod Gubarev worked on the translation in Vienna, supported
by the FWF, grant P28079. 
We thank Yurii A. Neretin for introducing us to the orbituary of I.A. Ado in Russian.

\end{document}